\documentclass[a4paper,11pt]{article}
\usepackage[psamsfonts]{amssymb}
\usepackage[psamsfonts]{eucal}

\usepackage{amsmath}
\usepackage{theorem}
\newtheorem{thm}{Theorem}
\newtheorem{prop}[thm]{Proposition}

\newtheorem{lm}[thm]{Lemma}

\newtheorem{rem}{Remark}

\theorembodyfont{\rmfamily}

\def\ds{\displaystyle}

\theorembodyfont{\rmfamily}
\begin{document}

\title{Compact complex surfaces with geometric structures related to split quaternions}

\author{Johann Davidov, Gueo Grantcharov, Oleg
Mushkarov, Miroslav Yotov}

\date{}
\maketitle

 \rm
\begin{abstract}

\smallskip

We study the problem of existence of geometric structures on
compact complex surfaces that are related to split quaternions.
These structures, called para-hypercomplex, para-hyperhermitian
and para-hyperk\"ahler, are analogs of the hypercomplex,
hyperhermitian and hyperk\"ahler structures in the definite case.
We show that a compact 4-manifold carries a para-hyperk\"ahler
structure iff it has a metric of split signature together  with
two parallel, null, orthogonal, pointwise linearly independent
vector fields. Every compact complex surface admitting a
para-hyperhermitian structure has vanishing first Chern class and
we show that, unlike the definite case, many of these surfaces
carry infinite dimensional families of such structures. We provide
also compact examples of complex surfaces with para-hyperhermitian
structures which are not locally conformally para-hyperk\"ahler.
Finally, we discuss the problem of non-existence of
para-hyperhermitian structures on Inoue surfaces of type $S^0$ and
provide a list of compact complex surfaces which could carry
para-hypercomplex structures.

\end{abstract}

\parindent 0cm
\vskip .1in

\section{Introduction}
    It was noticed long ago \cite{Ward} that many integrable systems
arise as reductions of self-dual Yang-Mills equations in signature $(2,2)$ and it is known that
they allow approaches via Lax pairs and twistor theory (see \cite{Dun-West} for a recent survey).
The geometry of a superstring with N=2 supersymmetry was shown in \cite{OV}, \cite{OV2} to be
described by a space-time with a pseudo-K\"ahler metric of signature (2,2), whose curvature
satisfies the (anti) self-duality equations. As noticed in \cite{OV}
this space-time also admits a (local) holomorphic (2,0)-form, parallel with respect to the Levi-Civita
connection. The structures obtained in this way define a holonomy reduction to the group
$SU(1,1)\cong SL(2,\mathbb{R})$ and are an indefinite analog of hyperk\"ahler structures which have
holonomy $SU(2)\cong Sp(1)$. Mathematically, these structures are described by quadruples
$(g,I,S,T)$ where $g$ is a signature $(2,2)$ metric and $I,S,T$ are parallel endomorphisms of the
tangent bundle such that:
 \begin{equation}\label{one}
 I^2=-S^2=-1,\quad T=IS=-SI, \quad g(IX,IY)=-g(SX,SY)= g(X,Y)
  \end{equation}
\smallskip
In the literature such structures are called hypersymplectic \cite{Hitchin}, neutral hyperk\"ahler
\cite{Kamada99}, para-hyperk\"ahler \cite{BV, DGMY}, pseudo-hyperk\"ahler \cite{Dun-West}, etc.
They are not preserved by a conformal change of the metric and a natural conformally invariant
generalization is to relax the condition for covariant constancy of $I,S,T$ to their integrability
(see Section 2). Such structures are an indefinite analog of the hyperhermitian structures and are
called para-hyperhermitian \cite{DGMY} or neutral hyperhermitian \cite{Kamada99, Kamada02}.  In
dimension 4, they are self-dual and, similarly to the positive definite case, always admit
connections with skew-symmetric torsion and holonomy $SL(2,\mathbb{R})$ \cite{kuche}. This geometry
was considered in \cite{Hull}, where it was argued that the (2,2) supersymmetric string based on
chiral multiplets is a theory of self-dual gravity. If we forget the metric $g$ and consider only
the triples $(I,S,T)$ of integrable endomorphisms of the tangent bundle satisfying the algebraic
conditions in (\ref{one}), the structures are called para-hypercomplex \cite{BV, DGMY}, neutral
hypercomplex \cite{Kamada99, Kamada02} or complex product \cite{Andrada, AnSal}. They provide
examples of geometric structures with special holonomy of a non-metric connection.

  \vspace{.1in}

Due to the non-elliptic nature of the self-duality equations in
the split signature case, their solutions are more flexible. For
example, a conformal self-dual or anti-self-dual structure is not
necessary analytic unlike the definite case. As a consequence,
most of the research deals with local properties of the
structures. However, topological information like the Kodaira
classification of compact complex surfaces allows one to study
global properties. Important examples in this direction are the
classifications of compact pseudo-K\"ahler Einstein and
para-hyperk\"ahler surfaces obtained
 by Petean \cite{P} and  Kamada \cite{Kamada99,
Kamada02}, respectively.

 \vspace{.1in}

In this paper, we study the compact 4-manifolds admitting para-hyperhermitian or para-hypercomplex
structures and our first aim is to relate the existence of para-hyperk\"ahler structures to the
existence of  parallel null orthogonal vector fields. More precisely, in Section 3 we show that if
a compact 4-manifold with a $(2,2)$-signature metric admits  two parallel, null, orthogonal,
pointwise linearly independent vector fields, then it is a torus or a primary Kodaira surface and
we notice that these surfaces do admit such vector fields (Theorem \ref{T}) .

\vspace{.1in}

A drastic difference between the definite and the split signature case is that some compact
para-hypercomplex $4$-manifolds do not admit compatible $(2,2)$-signature metrics, unlike the usual
hypercomplex manifolds \cite{DGMY}. We showed however that every compact para-hypercomplex
4-manifold (para-hypercomplex surface) has a double cover which admits a compatible
para-hyperhermitian metric \cite{DGMY}. Heuristically, this is due to the fact that
$GL(1,\mathbb{H}')/SU(1,1)= \mathbb{R}\backslash \{0\}$ has two connected components, where
$\mathbb{H}'$ is the algebra of split quaternions. Using the fact that the canonical bundle of a
complex surface with a para-hyperhermitian structure has a nowhere-vanishing smooth section, we
list in Section 4 the possible candidates for para-hyperhermitian surfaces (Theorem \ref{Ch0}). A
main observation in Theorem \ref{exist} is that most of these surfaces do admit para-hyperhermitian
structures which come in infinite dimensional families. This shows that the para-hyperhermitian
structures are much more flexible than the hyperhermitian ones.

  \vspace{.1in}

The considerations in Section 5 are motivated by the fact that, unlike the positive definite case,
there are compact para-hyperhermitian surfaces which are not locally conformally
para-hyperk\"ahler.  In Theorem \ref {lcpk}  we obtain a descrition of compact complex surfaces
admitting locally conformally para-hyperk\"ahler structures. To do this we first reduce the list of
possible candidates to those considered in Theorem \ref {exist} and then notice that the structures
constructed there, are in fact locally conformally para-hyperk\"ahler. An additional restriction
comes from the observation that the canonical bundle of such a surface is of real type in the sense
of \cite{ADl}. Moreover we give a construction leading to an infinite dimensional family of
para-hyperhermitian structures which are not locally conformally para-hyperk\"ahler.

  \vspace{.1in}

In Section 6 we provide a list of possible compact
para-hypercomplex surfaces by using Theorem \ref{Ch0} and the fact
that up to a double cover every para-hypercomplex surface is
para-hyperhermitian. Moreover, we construct a para-hypercomplex
structure on a surface in this list which does not admit a
compatible para-hyperhermitian metric.

  \vspace{.1in}

Finally, in Section 7 we study the Inoue surfaces of type $S^0$
which, as is well-known \cite{Hasegawa}, are solvmanifolds. We prove
in Theorem \ref{Inoue} that they do not admit para-hyperhermitian
structures with left-invariant canonical (2,0)-forms. This is a
slight generalization of the well-known result \cite{BV} that these
surfaces have no para-hyperhermitian structures induced by
left-invariant ones. This observation makes reasonable the
conjecture that the Inoue surfaces of type $S^0$ do not admit
para-hyperhermitian structures at all.

\vspace{.1in}

{\bf Acknowledgements:} We would like to thank M. Dunajski for his interest in the paper and are
especially grateful to V.Apostolov and G. Dloussky for several helpful discussions on problems
related to this project. Part of it was done during the visit of the first and the third-named
authors in the Abdus Salam School of Mathematical Sciences, GC University Lahore, Pakistan and the
second named author's visit to the Institute of Mathematics and Informatics at the Bulgarian
Academy of Sciences. The authors thank the two institutions for their hospitality.

\section{Preliminaries}
 Denote by $\mathbb{H}'$ the algebra of split quaternions, i.e.
 $$\mathbb{H}'=\{q=a+bi+cs+dt \in \mathbb{R}^4|~i^2=-1, s^2=t^2=1,
t=is=-si\} .$$ They are associated with a natural scalar product of
split signature (2,2) such that $|q|^2=a^2+b^2-c^2-d^2$. Based on
the algebra $\mathbb{H}'$, one defines an {\it almost para-hypercomplex}
structure on a manifold $M$ as a triple $(I,S,T)$ of
anti-commuting endomorphisms of the tangent bundle $TM$ with
$I^2=-Id$ and $S^2=T^2=Id$, $T=IS$. Such a structure is called
{\it para-hypercomplex} if $I,S,T$ satisfy the integrability condition
$N_I=N_S=N_T=0$, where
$$N_A(X,Y)=A^2[AX,AY]+[X,Y]-A[AX,Y]-A[X ,AY]$$ is the Nijenhuis
tensor associated with $A=I,S,T$. The para-hypercomplex structures
are the "split analog" of hypercomplex structures.

An almost product structure $S$ is integrable if and only if  the
eigenbundles $T^{\pm}=\{X\in TM: SX=\pm X\}$ are involutive
\cite{Walker}. Therefore, if $(I,S,T)$ is a para-hypercomplex
structure, then $T^{\pm}$ are two transversal involutive
distributions mapped to each other by the complex structure $I$.
Conversely, if a complex manifold $(M,I)$ admits such distributions
$T^{\pm}$, then we can define a para-hypercomplex structure setting
$S=Id$ on $T^{+}$, $S=-Id$ on $T^{-}$ and $T=IS$.

\medskip

A pseudo-Riemannian metric $g$ for which the endomorphisms $I,S,T$
are skew-symmetric is called {\it para-hyperhermitian}. Such a
metric necessarily has split signature and is also called neutral
hyperhermitian. Every para-hypercomplex structure on a 4-manifold
locally admits a para-hyperhermitian metric but a globally defined
one may not exist. More precisely, the following proposition is
true \cite{DGMY}.

\begin{prop}\label{cover}
Every para-hypercomplex structure on a 4-manifold $M$ determines a
conformal class of para-hyperhermitian  metrics up to a double cover
of $M$.
\end{prop}

Examples of para-hypercomplex structures that do not admit
para-hyperhermitian metrics are given in \cite{DGMY}. We provide
another example in Section 6.

\medskip

Para-hypercomplex four-manifolds can also be characterized in the following way:
\begin{prop} A four-dimensional smooth manifold  admits a para-hypercomplex structure if and only
if it admits two complex structures $I_1$ and $I_2$ yielding  the same
orientation and such that $I_1I_2+I_2I_1=2pId$ for a constant $p$ with $|p|>1$.
\end{prop}
{\bf Proof}. Suppose $I_1$ and $I_2$ are two complex structures such that  $I_1I_2+I_2I_1=2pId$
with $|p|>1$. Then $$I=I_1, ~S=\displaystyle{\frac{1}{2\sqrt{p^2-1}}}[I_1,
 I_2],~ T=-\displaystyle{\frac{1}{\sqrt{p^2-1}} }(I_1 +pI_2)$$
form an almost para-hypercomplex structure \cite{GL}. The
integrability of the structures $I$, $S$, $T$ is proved in
\cite[Lemma 1]{DGMY} based on the fact \cite{DGMY}  that, for each
point, there is a locally defined  para-hyperhermitian metric $g$.

  Conversely, suppose that we a given are para-hypercomplex structure $(I,S,T)$.
 Take any real number $p$ with $|p|>1$ and set
$$
I_1=I,\quad I_2=-pI-\sqrt{p^2-1}\,T.
$$
Then $[I_1,T]=2\sqrt{p^2-1}S$ and  $I_1I_2+I_2I_1=2pId$. It is
well-known that there is a unique torsion-free connection $\nabla$
such that $\nabla I=\nabla S=\nabla T=0$ (an analog of the Obata
connection) \cite[Theorem 3.1]{Andrada}. Clearly, $\nabla I_1=\nabla
I_2=0$. Since $\nabla$ is torsion-free, this implies that the
Nijenhuis tensors of $I_1$ and $I_2$ vanish, thus $I_1$ and $I_2$
are integrable. Take a point $x$ of the manifold and a metric on the
tangent space at $x$ which is compatible with $(I,S,T)$. Let $E$
be a non-isotropic tangent vector. Then $E_1=E$, $E_2=IE$, $E_3=SE$,
$E_4=TE$ is an orthogonal basis with $E_2=I_1E_1$, $E_4=I_1E_3$.
Moreover
$$I_2E_1=-pE_2-\sqrt{p^2-1}E_4, I_2E_3=-\sqrt{p^2-1}E_2-pE_4.$$ Therefore $I_1$ and $I_2$ determine
the same orientation. {\it Q.E.D.}

\begin{rem} If $I_1I_2+I_2I_1=2pId$ for a constant $p$ with $|p|<1$, then $I_1,I_2$ determine a usual hypercomplex structure and vice versa.
\end{rem}

\smallskip

Let $(g,I,S,T)$ be an almost para-hyperhermitian structure on a $4$-manifold. Then we can define
three fundamental $2$-forms $\Omega_i$, $i=1,2,3$, setting $$\Omega_1(X,Y)=g(IX,Y), \hspace{.1in}
\Omega_2(X,Y)=g(SX,Y), \hspace{.1in} \Omega_3(X,Y)=g(TX,Y).$$ Note that the form
$\Omega=\Omega_2+i\Omega_3$ is of type $(2,0)$ with respect to $I$. As in the definite case, the
corresponding Lee forms are defined by
$$\theta_1=\delta\Omega_1\circ I, \hspace{.1in}\theta_2=\delta\Omega_2\circ S,
\hspace{.1in}\theta_3=\delta\Omega_3\circ T,$$ where $\delta$ is the codifferential with respect to
$g$. It is well-known \cite{Boyer, GT, Iv-Zam, Kamada02} that $I,S,T$ are intergrable if and only
if $\theta_1=\theta_2=\theta_3$. Thus, for a para-hyperhermitian structure, we have just one Lee
form $\theta$; it satisfies the identities $$d\Omega_i=\theta\wedge\Omega_i,~ i=1,2,3.$$ When
additionally the three 2-forms $\Omega_i$ are closed, i.e. $\theta=0$, the para-hyperhermitian
structure is called {\it para-hyperk\"ahler} ( also hypersympectic or neutral hyperk\"ahler). When
$d\theta=0$ the structure is called {\it locally conformally para-hyperk\"ahler}. We note that, in
dimension 4,  the para-hyperhermitian metrics are self-dual and the para-hyperk\"ahler metrics are
self-dual and Ricci-flat \cite{Kamada02}. It is well-known that every hyperhermitian structure on a
4-dimensional compact manifold is locally conformally hyperk\"ahler \cite{Boyer}, but we shall see
in Theorem 9, that this is not true in the indefinite case.

A para-hyperhermitian $4$-manifold can be characterized by means
of the forms $\Omega_i$ and $\theta$ in the following way
\cite{Hitchin, Kamada02}.
\begin{prop}\label{phe}
Every para-hyperhermitian structure on a $4$-manifold is uniquely determined by three
non-degenerate $2$-forms $(\Omega_1, \Omega_2, \Omega_3)$  and a $1$-form $\theta$ such that
$$
-\Omega_1^2 = \Omega_2^2 = \Omega_3^2, \quad \Omega_l \wedge\Omega_m = 0, \, 1\leq l \neq m \leq 3,
\quad d\Omega_l=\theta\wedge\Omega_l.
$$
\end{prop}

\begin{prop}
Let $(M,J)$ be a simply connected complex surface that  carries a para-hyperhermitian structure $\{g,I,S,T\}$
with $I=J$.  If $\theta$ is the Lie form of this structure, then $dd^c\theta=0$ and the class of
$\theta$ in $$\frac{Ker(dd^c)}{Im(d)+Im(d^c)}$$ depends only on $J$.
\end{prop}
{\bf Proof}. The form  $\Omega=\Omega_2+i\Omega_3$ is of type
$(2,0)$ with respect to $I=J$ and nowhere-vanishing. Moreover
$d\Omega=\theta\wedge\Omega$, hence $d\theta\wedge\Omega=0$ which
implies $(d\theta)^{(0,2)}=0$. Then $(d\theta)^{(2,0)}=0$ since
$\theta$ is real-valued. Therefore $d\theta$ is a $(1,1)$-form. It
follows that  $\partial\theta^{(1,0)}=0$ and
$\overline{\partial}\theta^{(0,1)}=0$.  Then $dd^c\theta =
-i\partial\overline{\partial}\theta =
-i\partial\overline{\partial}\theta^{(1,0)} =
i\overline{\partial}(\partial\theta^{(1,0)})=0$.

  Let now $\{g',I',S',T'\}$ be another para-hyperhermitian structure on $(M,J)$ with $I'=J$.
Denote by $\Omega_1', \Omega_2',\Omega_3'$ the $2$-forms
determined by this structure and let $\theta'$ be the
corresponding Lee form. The  form $\Omega'=\Omega_2'+i\Omega_3'$
is of type $(2,0)$ with respect to $J$, hence $\Omega'= F\Omega$
for a nowhere-vanishing complex-valued smooth function $F$. Since
$M$ is simply connected, there is a smooth function $\varphi$ such
that $F=|F|e^{i\varphi}$. Then $(\theta'-\theta -d\ln
|F|-id\varphi)\wedge \Omega=0$ which implies $(\theta'-\theta
-d\ln |F|-id\varphi)^{(0,1)}=0$. It follows that $\theta'=\theta
+d\ln |F|+d^c\varphi$, so $\theta$ and $\theta'$ determine the
same class in
$$\frac{Ker(dd^c)}{Im(d)+Im(d^c)}.$$ {\it Q.E.D.}

\smallskip

Note that the cohomology class in Proposition 4 is related to the Aeppli cohomology groups (see
\cite{A}):
$$H^{p,q}_A=\frac{Ker(dd^c)\cap \Omega^{p,q}}{(Im(d)+Im(d^c))\cap\Omega^{p,q}}.$$

\smallskip

We can say a little bit more  for locally conformally
para-hyperk\"ahler structures. The first Chern class of a
holomorphic line bundle is determined by the coboundary map (the
Bockstein map)
$$\delta: H^1(M, \mathcal{O}^*) \rightarrow H^2(M, \mathbb{Z}),$$ where $\mathcal{O}^*$ is the
sheaf of non-vanishing holomorphic functions. The equivalence classes of topologically trivial
holomorphic line bundles are in the kernel $H^1_0(M, \mathcal{O}^*)$ of the map $\delta$.
 Then  following \cite{ADl} we consider the sequence of natural morphisms
$$ H^1(M,\mathbb{R})\mapsto H^1(M,\mathbb{R}_+)\hookrightarrow H^1(M,\mathbb{C}^*)\mapsto H^1_0(M,\mathcal{O}^*),$$
where the first morphism is induced by the exponential map $\mathbb{R}\mapsto \mathbb{R}^+$ and we
say that a bundle $\mathcal{L}\in H^1_0(M,\mathcal{O}^*) $ is of {\it real type} if its class is in
the image of $H^1(M,\mathbb{R}^+)$.

\begin{lm}\label{real}
If $M$ carries a locally conformally para-hyperk\"ahler structure,
then its canonical bundle is of real type.
\end{lm}
{\bf Proof}. Let $\Omega=\Omega_2+i\Omega_3$ be the  $(2,0)$-form  and $\theta$ be the  Lie form of
the given structure. Cover $M$ by open sets $\{U_{\alpha}\}$ such that the intersections
$U_{\alpha}\cap U_{\beta}$ are connected and $\theta|U_{\alpha}=d\varphi_{\alpha}$ for a smooth
function $\varphi_{\alpha}$. Then $\varphi_{\alpha}=\varphi_{\beta} +c_{\alpha\beta}$ on
$U_{\alpha}\cap U_{\beta}$, where $c_{\alpha\beta}$ are constants. We have
$d(e^{-\varphi_{\alpha}}\Omega)=0$, so $e^{-\varphi_{\alpha}}\Omega$ are local holomorphic sections
of the canonical bundle. These sections determine the transition functions
$\psi_{\alpha\beta}=e^{c_{\alpha\beta}}$. {\it Q.E.D.}

\section{Para-hyperk\"ahler surfaces and parallel null vector fields}

It has been shown by H. Kamada \cite{Kamada99, Kamada02} that the
only compact complex surfaces admitting para-hyperk\"ahler
structures are the primary Kodaira surfaces and the complex tori.
Moreover, he has described all such structures on these surfaces
in terms of the solutions of non-linear PDE's for a scalar
function \cite{Kamada99, Kamada02}. The aim of this section is to
find another characterization of para-hyperk\"ahler surfaces by
showing that they coincide with the compact 4-manifolds admitting
metrics of signature $(2,2)$ and pairs of parallel and orthogonal
null vector fields. Before stating our main result in this
direction we shall prove an auxiliary lemma.

\begin{lm}\label{prJ}
Let $M$ be a 4-manifold with a metric $g$ of signature $(2,2)$ and
let $X$ and $Y$ be orthogonal null vector fields which are
linearly independent at every point of $M$. Then the triple
$(g,X,Y)$ determines an orientation and a unique $g$- and orientation compatible almost complex
structure $J$ on $M$ such that $JX=Y$.
\end{lm}

{\bf Proof}. We first show that in a neighbourhood of every point of
$M$, there exist vector fields $Z,T$ such that:
\begin{enumerate}
\item[$(i)$] $(X,Y,Z,T)$ is a local frame of the tangent bundle $TM$;

\item[$(ii)$] $g(X,Z)=1$, $g(X,T)=0$;\quad $g(Y,Z)=0$, $g(Y,T)=1$.
\end{enumerate}

Indeed, by the Witt theorem, for every $p\in M$, there exist
isotropic tangent vectors $u,v\in T_pM$ such that $(X_p,Y_p,u,v)$
is a basis of $T_pM$ and $g(X_p,u)=1$, $g(X_p,v)=0$, $g(Y_p,u)=0$,
$g(Y_p,v)=1$. Extend $u,v$ to vector fields $U,V$ in a
neighbourhood of $p$ and consider the system

$$1=\alpha g(U,X) + \beta g(V,X),\quad 0=\alpha g(U,Y) + \beta
g(V,Y)$$

 with respect to the unknown functions $\alpha,\beta$. The
determinant of this system at the point $p$ is equal to $1$, hence
in a neighbourhood of $p$ it has a (unique) solution of smooth
functions $\alpha,\beta$. Similarly for the system $$ 0=\phi g(U,X)
+ \psi g(V,X),\quad 1=\phi g(U,Y) + \psi g(V,Y).$$
Set
$$
Z=\alpha U+\beta V,\quad T=\phi U +\psi V.
$$
We have $\alpha(p)=1$, $\beta(p)=0$, thus $Z_p=u$; similarly
$T_p=v$. Therefore $X,Y,Z,T$ form a frame of vector fields in a
neighbourhood of $p$.

Now let $\widetilde Z$, $\widetilde T$ be another pair of vector
fields around $p$ having the properties $(i)$ and $(ii)$ stated
above. Then they have the form
$$
\widetilde Z=a\,X+b\,Y +Z,\quad \widetilde T=c\,X+d\,Y+T,
$$
where $a,b,c,d$ are smooth functions. It follows that the frames
$(X,Y,Z,T)$ and $(X,Y,\widetilde Z,\widetilde T)$ determine the
same orientation. Thus, the orientation determined by $(X,Y,Z,T)$
does not depend on the choice of the vector fields $Z,T$ and we
shall say that it is determined by the triple $(g,X,Y)$.

 Next following \cite{DRRMML}, set $a=g(Z,Z)$, $b=g(T,T)$, $c=g(Z,T)$ and

\begin{equation}\label{ob}
\begin{array}{lll}
E_1=\displaystyle{\frac{1-a}{2}}X + Z, \quad E_2=\displaystyle{\frac{1-b}{2}}Y +T-cX, \\[8pt]
E_3=-\displaystyle{\frac{1+a}{2}}X + Z, \quad E_4=-\displaystyle{\frac{1+b}{2}}Y+ T - c X.
\end{array}
\end{equation}
Then $(E_1,E_2,E_3,E_4)$ is an orthogonal frame, positively
oriented with respect to the orientation determined by $(g,X,Y)$
and such that  $  g(E_1,E_1)=g(E_2,E_2)=1$,
$g(E_3,E_3)=g(E_4,E_4)=-1$. The almost complex structure $J$ for
which $JE_1=E_2$, $JE_3=E_4$ has the required properties. Let $K$
be another complex structure on $T_pM$ with these properties.
Define endomorphisms $S$ and $T$ of $T_pM$ such that $S^2=T^2=Id$
and $SE_1=E_3$, $SE_2=-E_4$, $TE_1=E_4$, $TE_2=E_3$. Set $I=J$.
Then $IS=-SI=T$. Since $K$ is compatible with the metric and
orientation it can be written as
$K=\lambda_1I+\lambda_2S+\lambda_3T$, where
$\lambda_1,\lambda_2,\lambda_3$ are real numbers with $
\lambda_1^2-\lambda_2^2-\lambda_3^2=1$. Moreover, in view of
(\ref{ob}), the identity $KX=Y$ is equivalent to
$KE_1-KE_3=E_2-E_4$ which implies that $\lambda_1=1$,
$\lambda_2=\lambda_3=0$. Thus $K=J$ which proves the lemma. {\it
Q.E.D.}

A complex structure compatible with a split signature metric and
preserving a null distribution of dimension 2 is called {\it proper}
for the null distribution \cite{M}. Lemma~\ref{prJ} shows that there
is a unique proper complex structure if we fix two orthogonal vector
fields of the distribution.

\begin{thm}\label{T}
Let $(M,g)$ be a compact $4$-manifold with a metric $g$ of
signature $(2,2)$. Suppose that $M$ admits two parallel and
orthogonal null vector fields $X$,$Y$, linearly independent at
every point of $M$ and let $J$ be the almost complex structure
determined by  $(g,X,Y)$ as in Lemma \ref{prJ}. Then:
\begin{itemize}
\item [$(i)$] The structure
$(g,J)$ is (pseudo) K\"ahler.
\item [$(ii)$] The metric $g$ is Ricci-flat.
\item [$(iii)$] $(M,J)$ is either a torus or a primary Kodaira surface.
\item [$(iv)$] $M$ admits a para-hyperk\"ahler structure with metric $g$ and complex structure
$I=J$.
\end{itemize}

 Conversely, every torus and every primary Kodaira surface $(M,J)$ admits a metric $g$ of signature $(2,2)$ and vector fields
$X$, $Y=JX$ which are parallel, orthogonal, null  and linearly
independent at every point of $M$.
\end{thm}
{\bf Proof}. Since $X$ and $Y$ are parallel, the proof of
\cite[Theorem 3]{Walker} shows that, around every point of $M$,
there are local coordinates $(x,y,z,t)$, such that $\ds
X=\frac{\partial}{\partial x}$, $\ds Y=\frac{\partial}{\partial y}$
and the metric $g$ in these coordinates has the form

\begin{equation}\label{Walker metric}
g_{(x,y,z,t)} =\left(\begin{array}{cccc}
    0&0&1&0\\
    0&0&0&1\\
    1&0&a&c\\
    0&1&c&b\\
\end{array}\right)
\end{equation}
where $a$, $b$, $c$ are smooth functions independent of the
coordinates $x$ and $y$. According to Lemma~\ref{prJ}, the manifold
$M$ admits a unique almost complex structure $J$ compatible with $g$
and such that $\ds J\frac{\partial}{\partial
x}=\frac{\partial}{\partial y}$. Since $a$, $b$, $c$ do not depend
on $x$, $y$, it follows from \cite[Corollary 14]{DRRMML} that the structure $(g,J)$ is (pseudo) K\"ahler and
Ricci-flat . Then, by \cite[Corollary
2]{P}, $M$ is one of the following: a torus, a primary Kodaira
surface or a hyperelliptic surface. To prove the result we need first to
exclude the last option. Note that $X-iY$ is a parallel,
isotropic, nowhere-vanishing, holomorphic vector field on $M$. Every
hyperelliptic surface $M$ is the quotient of a product $E\times F$
of two elliptic curves by a finite fixed-point-free abelian group of
automorphisms. We can take $E$ of the form $\mathbb{C}/\Lambda$
where $\Lambda$ is the lattice generated by $1$ and a complex number
$\tau$ with $Im\,\tau>0$. Then $M$ is the quotient of $E\times F$ by
the group generated by certain translations of $F$ and a map of the
type $$\varphi(z,w)=(z+\tau/m, e^{2k\pi i/m}w),$$ where $m$ is one
of the numbers $2,3,4,6$ and $k\in\{1,...,m-1\}$
\cite[VI.19, VI.20]{B}. The quotient map $E\times F\to M$ is a
(finite) covering and we take the pull-back $g'$ of the metric $g$,
then lift $g'$ to the universal covering $\mathbb{C}^2$ of $E\times
F$. In this way we get a (pseudo) K\"ahler, Ricci-flat metric
$\widetilde g$ on $\mathbb{C}^2$. It is of the form
$$\widetilde g=\alpha\, dzd\overline{z} + 2Re(\gamma\,
dzd\overline{w}) + \beta\, dwd\overline{w}$$ for real smooth
functions $\alpha$, $\beta$ and a complex smooth function
$\gamma$. The lift $\ds\widetilde
U=\lambda\frac{\partial}{\partial z}+\mu\frac{\partial}{\partial
w}$ on $\mathbb{C}^2$ of a holomorphic vector field $U$ on $M$
satisfies the identity $\varphi_{\ast}\circ\widetilde U=\widetilde
U\circ\varphi$. This implies $\mu=0$ since
$\ds\varphi_{\ast}(\frac{\partial}{\partial
z})=\frac{\partial}{\partial z}$ and
$\ds\varphi_{\ast}(\frac{\partial}{\partial w})=e^{2k\pi
i/m}\frac{\partial}{\partial w}$. Therefore, if $Z$ is the
holomorphic vector field on $M$ given in the local coordinates
$(z,w)$ as $\displaystyle{\frac{\partial}{\partial z}}$, we have
$U=fZ$ for a function $f$. The function $f$ is holomorphic on the
compact manifold $M$, hence it is a constant. Thus, every
holomorphic vector field on $M$ is proportional to $Z$. It follows
that the vector field $\ds\partial_z=\frac{\partial}{\partial z}$
on $\mathbb{C}^2$ is parallel and null with respect to $\widetilde
g$. The fact that this field is null implies $\alpha=0$. Since it
is parallel, the Lie derivative of the  K\"ahler form
$$\Omega=-i(\alpha\, dz\wedge d\overline{z} + \gamma\, dz\wedge
d\overline{w}+ \overline{\gamma}\,dw\wedge d\overline{z}+ \beta\, dw\wedge d\overline{w})$$
vanishes. Thus, it follows by Cartan's formula that $$0= {\mathcal L}_{\partial_z}
\Omega=d\circ\imath_{\partial_z}\Omega=-id(\gamma d\overline{w})$$ since $d\Omega=0$ and
$\alpha=0$. This implies that the derivatives of $\gamma$ with respect to $z$, $\overline{z}$, $w$
vanish. Thus $\gamma$ depends only on the variable $w$ and $\gamma(w)$ is an anti-holomorphic
function. Then, since $\widetilde g$ is the lift of a metric on $E\times F$, $\gamma$ descends to
an anti-holomorphic function $\gamma'$ on $F$. By the maximum principle, $\gamma'\equiv const$,
therefore $\gamma$ is a constant. This constant is not zero since $\widetilde g$ is non-degenerate.
Then $\widetilde g$ is not invariant under the map $\varphi$, so it does not descend to $M$, a
contradiction.

Next we show that the primary Kodaira surfaces and 4-tori do admit
para-hyperk\"ahler structures, compatible with $g$ and $I=J$.
Suppose that $M$ is a primary Kodaira surface. Then it can be
obtained in the following way. Consider the affine transformations
$\varphi_k(z,w)$ of ${\Bbb C}^2$ given by
$$\varphi_k(z,w) = (z+a_k,w+\overline{a}_kz+b_k),$$ where $a_k$, $b _k$, $k=1,2,3,4$, are complex
numbers such that $$a_1=a_2=0, \quad Im(a_3{\overline a}_4) =
b_1\neq 0,\quad b_2\neq 0.$$ They generate a group $G$ of affine
transformations acting freely and properly discontinuously on
${\Bbb C}^2$ and $M$ is the  quotient space ${\Bbb C}^2/G$ for a
suitable choice of $a_k$ and $b _k$ \cite[p.786]{Kodaira}. Taking
into account the identities
$$\varphi_{k\ast}(\frac{\partial}{\partial
z})=\frac{\partial}{\partial z} + \overline{a_k}
\frac{\partial}{\partial w},\quad
\varphi_{k\ast}(\frac{\partial}{\partial
w})=\frac{\partial}{\partial w}$$ we see that every holomorphic
vector field on $M$ is proportional to the vector field $W$ given
in the local coordinates $(z,w)$ as $\ds\frac{\partial}{\partial
w}$. Therefore the vector field
$\ds\partial_w=\frac{\partial}{\partial w}$ on $\mathbb{C}^2$ is
parallel and null with respect to the lift
$$\widetilde g=\alpha\, dzd\overline{z} + 2Re(\gamma\,
dzd\overline{w}) + \beta\, dwd\overline{w}$$ of the metric $g$.
Arguments similar to that above show that $\beta=0$ and $\gamma=
const\neq 0$. Then the K\"ahler form of the K\"ahler metric
$\widetilde g$ is given by
$$\widetilde\Omega_1= -i(\alpha(z)\, dz\wedge d\overline{z} + \gamma\, dz\wedge d\overline{w}+
\overline{\gamma}\,dw\wedge d\overline{z}),$$ where $\alpha(z)$ is a
smooth function depending only of $z$. Set $$\widetilde
\Omega_2=\gamma \, dz\wedge dw+\overline{\gamma}\,
d\overline{z}\wedge d\overline{z}, \quad \widetilde
\Omega_3=-i(\gamma \, dz\wedge dw-\overline{\gamma}\,
d\overline{z}\wedge d\overline{z}).$$ By Proposition~\ref{phe}, the
forms $\widetilde\Omega_1$, $\widetilde\Omega_2$,
$\widetilde\Omega_3$ determine a para-hyperk\"ahler structure on
$\mathbb{C}^2$. Since these $2$-forms are invariant under the action
of the group $G$, we obtain a para-hyperk\"ahler structure on $M$
with metric $g$ and complex structure identical to the complex
structure of $M$.

Now consider the case when $M$ is a complex torus. Let
${\partial_z}$ and ${\partial_w}$ be the global holomorphic vector
fields given in the standard local coordinates $(z,w)$ on the torus
by $\ds{\partial_z}=\frac{\partial}{\partial z}$,
$\ds{\partial_w}=\frac{\partial}{\partial w}$. Then the holomorphic
vector field $U=X-iY$ on $M$ is a linear combination $U=\lambda
{\partial_z}+\mu {\partial_w}$, where $\lambda$ and $\mu$ are
constants. The Lie derivative with respect to $U$ of the K\"ahler
form $$\Omega_1=-i(\alpha\, dz\wedge d\overline{z} + \gamma\,
dz\wedge d\overline{w}+ \overline{\gamma}\,dw\wedge d\overline{z}+
\beta\, dw\wedge d\overline{w})$$ of the metric $g$ vanishes. It
follows that the derivatives of $\lambda\alpha+\mu\overline{\gamma}$
with respect to $z$,$w$, $\overline w$ vanish and the derivatives of
$\lambda\gamma+\mu\beta$ with respect to $z$, $\overline z$, $w$
also vanish. Therefore the derivatives of the functions
$f^{\pm}=(\overline{\lambda}\alpha+\overline{\mu}\gamma)\pm(\overline{\lambda}\overline{\gamma}+\overline{\mu}\beta)$
with respect to $\overline{z}$ and $\overline{w}$ vanish. Then
$f^{\pm}dz\wedge dw$ are globally defined closed forms. We have
$\alpha
|\lambda|^2+\gamma\lambda\overline{\mu}+\overline{\gamma}\overline{\lambda}\mu+\beta
|\mu|^2=0$ since the vector field $U$ is isotropic. This identity
implies that
$|f^{\pm}|^2=|\lambda\mp\mu|^2(|\gamma|^2-\alpha\beta)$. Note that either
$\lambda+\mu\neq 0$ or $\lambda-\mu\neq 0$ since $U\neq 0$. In the
first case we define real-valued $2$-forms $\Omega_2$ and $\Omega_3$
by $$\Omega_2+i\Omega_3=\frac{2f^{-}}{|\lambda+\mu|}dz\wedge dw$$
and in the second case we set
$$\Omega_2+i\Omega_3=\frac{2f^{+}}{|\lambda-\mu|}dz\wedge dw.$$ In
both cases the forms $\Omega_1$, $\Omega_2$, $\Omega_3$ determine a
para-hyperk\"ahler structure on $M$ with metric $g$ and complex
structure $I=J$.

Finally we show that the primary Kodaira surfaces and 4-tori admit
metrics with 2 parallel orthogonal null vector fields.  Let $M$ be
a primary Kodaira surface represented as ${\Bbb C}^2/G$, where the
group $G$ has been described above. As in \cite{P}, set
$\alpha(z)=f(z) -\gamma\,z-\overline{\gamma}\,\overline{z}$, where
$f(z)$ is a smooth function on ${\Bbb C}$ satisfying the
identities $f(z+a_3)=f(z)$, $f(z+a_4)=f(z)$. Then the metric
$\widetilde g=\alpha\, dzd\overline{z} + 2Re(dzd\overline{w})$ on
$\mathbb{C}^2$ descends to a split signature K\"ahler, Ricci flat
metric $g$ on $M$ for which the holomorphic vector field $W$ is
parallel and null (the metric $g$ is flat if $f\equiv const$).
Hence the real and imaginary parts of the vector field $W$ have
the required properties.

The case of a complex torus $M$ is similar. Let $M$ be the
quotient of $\mathbb{C}^2$ by a lattice $<a_1,a_2,a_3,a_4>$. Take
a smooth function $\alpha$ on ${\Bbb C}$ such that
$\alpha(z+a_3)=\alpha(z)$, $\alpha(z+a_4)=\alpha(z)$. Then the
metric $$\widetilde g=\alpha\, dzd\overline{z} +
2Re(dzd\overline{w})$$ on $\mathbb{C}^2$ descends to a a split
signature  Ricci flat, K\"ahler  metric $g$ on the torus $M$
\cite{P}. For this metric, the holomorphic vector field
$\ds\frac{\partial}{\partial w}$ on $M$ is parallel and null. {\it
Q.E.D.}

\section{Para-hyperhermitian structures and compact complex surfaces with vanishing first Chern class}

Let $(g,I,S,T)$ be a para-hyperhermitian structure on a
compact $4$-manifold $M$ with fundamental $2$-forms $\Omega_i$,
$i=1,2,3$. Then $\Omega=\Omega_2+i\Omega_3$ is of type $(2,0)$
w.r.t. the complex structure  $I$, thus the canonical bundle of
the complex manifold $(M,I)$ is smoothly trivial. For a para-hyperk\"ahler
structure the $(2,0)$-form $\Omega$ is holomorphic, hence the
canonical bundle of $(M,I)$ is holomorphically trivial. Thus a
compact complex surface $(M,J)$ admits a para-hyperhermitian
structure with $I=J$ only if its integral first Chern class
$c_1(J)$ vanishes. A weaker condition is the vanishing of the real
first Chern class $c_1^{\mathbb{R}}(J)$. If $c_1^{\mathbb{R}}(J)=0$,
then $c_1(J)$ is torsion and the canonical bundle of $(M, J)$ is
flat as well as the principle circle bundle corresponding to it.
Since the flat principal $G$-bundles for any Lie group $G$ are in
bijection with the conjugacy classes of $G$-representations of the
fundamental group of $M$, there is a finite covering
$(\widetilde{M},\widetilde{J})$ of $(M,J)$ with
$c_1(\widetilde{J})=0$.

The classification of surfaces with topologically trivial canonical bundle seems to be known to the
experts, but we were not able to find an explicit proof. So we provide a short proof below for the
sake of completeness. The primary Kodaira surfaces from the list were defined in Theorem~\ref{T}
while the Hopf surfaces, the Inoue surfaces,  and the minimal properly elliptic surfaces of odd
first Betti number are discussed in Theorems~\ref{exist} and \ref{lcpk}. Note that we use the
notations for Inoue surfaces from \cite{Wall} and $S^o$ corresponds to $S_M$ and $S^{\pm}$
corresponds to $S^{\pm}_N$ in \cite{Inoue, ADl}.

\begin{thm}\label{Ch0}
Let $(M,J)$ be a compact complex surface with topologically trivial canonical bundle. Then $(M,J)$
is one of the following: a complex torus, a $K3$ surface, a primary Kodaira surface, a Hopf
surface, an Inoue surface of type $S^0$ or $S^{\pm}$ without curves, or a minimal properly elliptic
surface of odd first Betti number.
\end{thm}
{\bf Proof}. The proof is based on \cite{Wall} and makes use of the completed proof (c.f. \cite{Te}
and \cite{Li-Yau}) of a result due to Bogomolov \cite{Bog} about surfaces of class $VII$ with
vanishing second Betti number.

It is well-known that if the first Betti number $b_1$ of $M$ is
even, then it admits a K\"ahler metric.  In this case the
canonical bundle ${\cal K}$ of $M$ is not just topologically but
also holomorphically trivial. Indeed, by \cite[Th\'{e}or\`{e}me
1]{B1}, $M$ admits a finite holomorphic covering with
holomorphically trivial canonical bundle. Then a tensor power
${\cal K}^d$ of ${\cal K}$ is also holomorphically trivial (c.f.,
for example, \cite[(16.2)Lemma, p.54]{BHPV}). Since ${\cal K}$ is
topologicallay trivial, it follows that it is holomorphically
trivial.  Then $M$ is either a torus or a $K3$ surface by the
Kodaira-Enriques classification. So we further consider only the
case $b_1$ odd . Notice that, since $c_1(M)=0$, the adjunction
formula shows that $M$ is minimal. Hence it also satisfies
$c_1^2(M)=0$. Such surfaces are easy to identify in  the
Kodaira-Enriques classification. In particular their Kodaira
dimension $k$ is one the numbers $-\infty, 0,1$. We start with
$k=-\infty$. In this case it follows from \cite[Sec. 6]{Wall} that
$c_1^2=-b_2(M)$, so $b_2(M)=0$. Then $c_1^{\mathbb{R}}(M)=0$ and
by the above mentioned result of Bogomolov, $(M,J)$ is either a
Hopf surface or Inoue surface. The second case is $k=0$. Then we
see, again from the Kodaira-Enriques classification, that $(M,J)$
is a Kodaira surface. Lastly, if $k=1$, then $(M,J)$ is a properly
elliptic surface.  {\it Q.E.D.}

\smallskip

\begin{rem} It is well known that every torus, $K3$ surface,  primary
Kodaira surface and primary Hopf surface  has a topologically trivial canonical bundle \cite{BHPV}.

\smallskip

Some of the non-primary Hopf surfaces have trivial canonical bundle, but some of them do not.  For
example, the quotient of a quaternionic Hopf surface by a finite subgroup of $SU(2)$ admits a
hypercomplex structure, hence it has  a non-vanishing section of its canonical bundle, so the
latter is trivial. We can define a Hopf surface with topologically non-trivial canonical bundle as
follows. Let $(z,w)$ be the standard coordinates on ${\Bbb C}^2\setminus(0,0)$ and $G$ the group of
transformations of  ${\Bbb C}^2\setminus(0,0)$ generated by
$g_0(z,w)=\displaystyle{(\frac{1}{2}z,\frac{1}{2}w)}$,
$§g_1(z,w)=(w,z)$. Consider the secondary Hopf surface $M=({\Bbb C}^2\setminus(0,0))/G$. Let $\mu$
be the representation of the fundamental group $G$ of $M$ yielding its canonical bundle ${\cal K}$.
According to \cite[p. 271]{Inoue}, the integral first Chern class of ${\cal K}$ vanishes if and
only if $\mu|Tor(G/[G,G])=1$. It is clear that $g_0$ is the free generator of $G$, $g_1^2=1$,
$g_0g_1g_0^{-1}g_1^{-1}=Id$ and $g_1 mod([G,G])\in Tor(G/[G.G])$. We have $\mu(g_1)=-1$, hence
${\cal K}$ is not topologically trivial (but note that ${\cal K}^2$ is trivial).

\smallskip

The integral first Chern class of a properly elliptic surface depends on the invariant $c(\eta)$
defined  in \cite[p. 140]{Wall}. The first Chern class is proportional to a generator $\rm c$ of
the second cohomology of the base of the elliptic fibration. In particular it
vanishes iff $c(\eta)$ is primitive, i.e. $c(\eta)=\pm c$.

\smallskip

The Inoue surfaces also have topologically trivial canonical bundle - the surfaces of type $S^+$
and $S^0$ admit non-vanishing smooth (2,0)-forms -see, for example, the proofs of
Theorems~~\ref{exist} and \ref{lcpk}, while to see that the integral first Chern class of an Inoue
surface of type $S^-$ vanishes we can use the condition given in \cite[p. 271]{Inoue}. Every Inoue
surface $S$ of type $S^{-}$ is the quotient of ${\mathbb C}\times{\bf H}$ by a group $G$ of
transformations acting freely and properly discontinuous (see the proof of Theorem~\ref{exist} for
a detailed description). The canonical bundle of $S$ is the associated bundle $({\mathbb
C}\times{\bf H})\times_{\mu}{\mathbb C}$, where the representation $\mu:G\to End({\mathbb C})$ is
defined by $g^{\ast}(dz\wedge dw)=\mu(g)(dz\wedge dw)$, $g\in G$, $(z,w)$ being the standard
coordinates on ${\mathbb C}\times{\bf H}$.  By \cite[p. 279]{Inoue}, $Tor(G/[G,G])$ is generated by
the transformations $g_k(z,w)=(z+b_kw+c_k,w+a_k)$, $k=1,2,$,
$g_3(z,w)=(z+\displaystyle{\frac{b_1a_2-b_2a_1}{r}},w)$ where $a_k, b_k, c_k, r$ are certain
numbers. These transformations leave invariant the form $dz\wedge dw$, thus $\mu|Tor(G/[G,G])=1$.
Therfore the canonical bundle of $S$ is topologically trivial by  \cite[p. 271]{Inoue}.

\end{rem}

\smallskip

In the next theorem we show that most of the surfaces listed in
Theorem \ref {Ch0} do admit para-hyperhermitian structures.
Moreover, they vary in infinite dimensional families.

\begin{thm}\label{exist}
The following compact complex surfaces admit infinite dimensional
families of para-hyperhermitian structures: complex tori, primary
Kodaira surfaces,  Inoue surfaces of type $S^{+}$, a special type of
minimal properly elliptic surfaces with odd first Betti number and
quaternionic primary Hopf surfaces. The complex tori and the primary Kodaira surfaces admit an infinite dimensional family of
non-locally conformally para-hyperk\"ahler structures.
\end{thm}
{\bf Proof}. The proof is case by case.

\smallskip

(1) {\it Complex tori}.

\smallskip

As we have mentioned, every complex torus of dimension 2 admits a para-hyperk\"ahler structure
\cite{Kamada02}. Here we shall construct an infinite dimensional family of para-hyperhermitian
structures which are not locally conformally para-hyperk\"ahler.

\smallskip

Let $M=\mathbb{C}^2/\Lambda$ be a complex torus with lattice $\Lambda$ generated by vectors $\tau_1,...,\tau_4$ where
$\tau_1=(a_1,0)$ and $\tau_2=(a_2,0)$. Take a non-constant real-valued smooth doubly-periodic
function $\varphi$ on ${\Bbb C}$ with periods $a_1$ and $a_2$. Let $(z,w)$ be the
standard coordinates on $\mathbb{C}^2$. Set
$$
\Omega_1=Im\,(e^{i\varphi}dz\wedge d\overline{w}),\quad
\Omega_2+i\Omega_3=e^{i\varphi} dz\wedge dw,\quad
\theta=i\frac{\partial\varphi}{\partial\overline{z}}d\overline{z}
-i\frac{\partial\varphi}{\partial z}dz.
$$
These forms descend to $M$ and, in view of Proposition~\ref{phe},
they determine a para-hyperhermitian structure on $M$ which is not
para-hyperk\"ahler. We have $d\theta=0$ exactly when the function
$\varphi$ is harmonic. In this case $\varphi$ is constant since it
is bounded. Hence the para-hyperhermitian structure on $M$ defined
above is not locally conformally para-hyperk\"ahler.

\smallskip

(2) {\it Primary Kodaira surface}s.

\smallskip

Every such a surface admits a para-hyperk\"ahler structure \cite{Kamada99, Kamada02}.

\smallskip

Here we shall use the description of primary Kodaira surfaces as
quotients $M=\mathbb{C}^2/G$ given in the proof of Theorem~\ref{T}.
Note first that the complex numbers $a_3$ and $a_4$ are linearly
independent over $\mathbb{R}$ since $Im(a_3{\overline a}_4)\neq 0$.
Now take a non-constant real-valued doubly-periodic function
$\varphi$ on $\mathbb{C}$ with periods $a_3$ and $a_4$. Set
$$
\Omega_1 = Im(e^{i\varphi}dz\wedge
d\overline{w})+iRe(e^{i\varphi}z)dz\wedge d\overline{z},\quad
\Omega_2+i\Omega_3 = e^{i\varphi} dz\wedge dw,\quad
\theta=i\frac{\partial\varphi}{\partial\overline{z}}d\overline{z}
-i\frac{\partial\varphi}{\partial z}dz,
$$
where $(z,w)$ are the standard coordinates on $\mathbb{C}^2$. Then
these forms satisfy the identities of Proposition~\ref{phe}.
Moreover, the forms $\Omega_1,\Omega_2,\Omega_3,\theta$ are
invariant under the action of the group $G$, so they define a
para-hyperhermitian structure on $M$ which is not locally
conformally para-hyperk\"ahler.

\smallskip

(3) {\it Quaternionic Hopf surfaces}.

\smallskip

These surfaces are the quotient spaces $M=(\mathbb{H}'\backslash\{0\})/\mathbb{Z}$, the action of
$\mathbb{Z}$ being generated by $L_a:q\rightarrow aq$, where $a$ is a fixed complex number with
$|a|>1$. If $q=z_1+sz_2$,  the action is $(z_1,z_2)\rightarrow(az_1, \overline{a}z_2)$. Then the
following forms define conformally para-hyperk\"ahler structure on $\mathbb{H}'\backslash\{0\}$
which descends to a locally conformal para-hyperk\"ahler structure on $M$ :
$$
\Omega_1=i\frac{dz_1\wedge d\overline{z}_1-dz_2\wedge
d\overline{z}_2}{|z_1|^2+|z_2|^2},\quad
\Omega_2+i\Omega_3=\frac{dz_1\wedge dz_2}{|z_1|^2+|z_2|^2},
$$
$$
\theta=-\frac{1}{|z_1|^2+|z_2|^2}(\overline{z}_1dz_1+z_1d\overline{z}_1+\overline{z}_2dz_2+z_2d\overline{z}_2).
$$
We can also take
$$
\Omega_1= i\frac{dz_1\wedge d\overline{z}_1-dz_2\wedge
d\overline{z}_2}{|z_1|^2+|z_2|^2}+i\partial\overline{\partial}\varphi,
$$
with a smooth real-valued function $\varphi$ depending only on
$z_1$.

\smallskip

(4) {\it Inoue surfaces of type $S^{+}$}

\smallskip

We first recall the construction of the Inoue surfaces of type $S^{\pm}$ \cite{Inoue}. Set
$\varepsilon=\pm 1$ and take a matrix $N=(n_{ij})\in GL(2,{\Bbb Z})$ with $det\,N=\varepsilon$
having two real eigenvalues $\alpha>1$ and $\varepsilon\alpha$. Note that $\alpha$ is a irrational
number. Choose real eigenvectors $(a_1,a_2)$ and $(b_1,b_2)$ corresponding to $\alpha$ and
$\varepsilon\alpha$, respectively. Take integers $p,q,r$, $r\neq 0$ and a complex number $t$. Let
$(c_1,c_2)$ be the solution of the equation
\begin{equation}\label{cc}
\varepsilon
(c_1,c_2)=(c_1,c_2)N^{tr}+(e_1,e_2)+\frac{b_1a_2-b_2a_1}{r}(p,q)
\end{equation}
where $N^{tr}$ is the transpose matrix of $N$ and
$$
e_k=\frac{1}{2}n_{k1}(n_{k1}-1)a_1b_1+\frac{1}{2}n_{k2}(n_{k2}-1)a_2b_2+n_{k1}n_{k2}b_1a_2,
\quad k=1,2.
$$
Let $G^{\varepsilon}$ be the group generated by the following
automorphisms of ${\Bbb C}\times {\bf H}$, ${\bf H}$ being the
upper half-plane:
\begin{equation}\label{G}
\begin{array}{l}
g_0=(z,w)=(\varepsilon z+\frac{1}{2}(1+\varepsilon)t,\alpha w)\\[6pt]
g_k(z,w)=(z+b_kw+c_k,w+a_k),\, k=1,2,\quad
g_3(z,w)=(z+\displaystyle{\frac{b_1a_2-b_2a_1}{r}},w).
\end{array}
\end{equation}
The group $G^{\epsilon}$ acts properly discontinuously and without
fixed points in view of (\ref{cc}) and the fact that $(a_1,b_1)$ and
$(a_2,b_2)$ are linearly independent vectors. The quotient $({\Bbb
C}\times {\bf H})/G^{\varepsilon}$ is a compact complex surface,
known as an Inoue surface of type $S^{\varepsilon}$.

 Given an Inoue surface of type $S^{+}$, we set $t_2=Im\,t$ and
\begin{equation}\label{s+}
\alpha_1=dx-\frac{1}{v}(y-t_2\frac{\ln v}{\ln\alpha})du,\quad
\alpha_2=dy-\frac{1}{v}(y-t_2\frac{\ln v}{\ln\alpha})dv,\quad
\alpha_3=\frac{du}{v}, \quad \alpha_4=\frac{dv}{v},
\end{equation}
where $z=x+iy$ and $w=u+iv$. These forms are linearly independent
and invariant under the action of the group $G^{+}$. Moreover
$$
d\alpha_1=\alpha_3\wedge\alpha_2-\frac{t_2}{v^2\ln\alpha}du\wedge
dv,\quad d\alpha_2=\alpha_4\wedge\alpha_2,\quad
d\alpha_3=\alpha_3\wedge\alpha_4,\quad d\alpha_4=0.
$$
Set
$$
\Omega_1=\alpha_1\wedge\alpha_3+\alpha_2\wedge\alpha_4,\quad
\Omega_2=\alpha_1\wedge\alpha_3-\alpha_2\wedge\alpha_4,\quad
\Omega_3=\alpha_1\wedge\alpha_4+\alpha_2\wedge\alpha_3.
$$
Then
$$
-\Omega_1^2 = \Omega_2^2 =
\Omega_3^2=2\alpha_1\wedge\alpha_2\wedge\alpha_3\wedge\alpha_4,
\quad \Omega_l \wedge\Omega_m = 0, \,1\leq l \neq m\leq 3, \quad
d\Omega_l=-\alpha_4\wedge\Omega_l.
$$
Therefore, by Proposition~\ref{phe},  $\Omega_1,\Omega_2,\Omega_3$
define an $G^{+}$-invariant para-hyperhermitian structure on ${\Bbb
C}\times {\bf H}$ which is locally conformally para-hyperk\"ahler
since its Lie form $\theta=-\alpha_4$ is closed. This structure
descends to a para-hyperhermitian structure on the Inoue surface
$S^{+}$.

We can deform $\Omega_1$ to
$\Omega_1+\displaystyle{\frac{i\partial\overline{\partial}\varphi}{Im(w)}}$
for arbitrary function $\varphi$ depending only on $Im(w)$ and
satisfying $\varphi(\alpha x)=\varphi(x)$. These functions are in
one-to-one correspondence with the functions on the circle
$S^1=\mathbb{R}^+/<\alpha >$ and form an infinite dimensional
family.
\bigskip

(5) {\it Minimal properly elliptic surfaces of odd first Betti
number}

\smallskip

A properly elliptic surface is, by definition, a compact complex
surface admitting a fibration $\pi:M\rightarrow B$ onto a complex
orbifold curve of genus $g\geq 2$ and generic fiber an elliptic
curve. Every such a surface is of Kodaira dimension 1 and has
universal  cover $\mathbb{C}\times {\bf H}$. Among these surfaces,
the ones with vanishing first Chern class are precisely those with
odd first Betti number. In this case $M$ has no singular fibers
\cite[Lemma 7.2]{Wall} and has a good orbifold base - see the
considerations preceding Theorem 7.4 in \cite{Wall}.

It is convenient to use here the description of the minimal
properly elliptic surfaces $M$ with odd first Betti number given
in \cite{Maehara}. Set
 $$
 D=\{(x,y)\in {\mathbb C}^2|~
Im(x/y)>0\}.
$$
According to \cite[Theorem 1]{Maehara}, every minimal elliptic surface of odd first Betti number is
a quotient of this (non-simply connected) domain by a discrete group $\Gamma$ generated by a finite
number of linear transformations of ${\mathbb C}^2$ of the form $L=\lambda M$, where
$\lambda\in\mathbb{C}\setminus\{0\}$ and $M \in SL(2,{\mathbb R})$. The matrices $L$ satisfy a
number of relations and the elliptic fibration is determined by the map $\pi: (x,y)\rightarrow
\displaystyle{\frac{x}{y}}$ on the covering space $D$.

Take a transformation   $L(x,y) = (\lambda(ax+by),\lambda(cx+dy))$ of $D$ where $\left(
\begin{array}{ll}
a& b\\
c& d\\
\end{array}\right) \in SL(2,\mathbb{R})$ and $\lambda\in
\mathbb{C}\setminus\{0\}$. Then
$$
L^{\ast}(dx\wedge dy)=|\lambda|^2dx\wedge dy,\quad
Im(\lambda(ax+by)\overline{\lambda(cx+dy)})=|\lambda|^2Im(x\overline{y}).
$$
Thus, when the number $\lambda$ is real, or equivalently when the
group $\Gamma$ is a subgroup of $GL(2, {\mathbb R})$, the forms
$$
\Omega_1 = Im(\frac{dx\wedge
d\overline{y}}{Im(x\overline{y})}),\quad \Omega_2+i\Omega_3 =
\frac{dx\wedge dy}{Im(x\overline{y})},\quad
\theta=\frac{1}{2iIm(x\overline{y})}(yd\overline{x}-\overline{y}dx+\overline{x}dy-xd\overline{y}).
$$
descend to $M$. It is easy to check that they determine a locally conformally para-hyperk\"ahler
structure.

 Taking a smooth real-valued function $\varphi$ such that
$\varphi(|\lambda|^2z)=\varphi(z)$ for all $z\in{\bf H}$, we obtain
an infinite-dimensional family of para-hyperhermitian structures on
$M$ by changing $\Omega_1$ to $$\Omega_1=
\displaystyle{Im(\frac{dx\wedge d\overline{y}}{Im(x\overline{y})})
+dd^c\varphi}.$$ Note that every function $\varphi$ depending only
on $arg(z)$ satisfies the above condition. {\it Q.E.D.}

\smallskip

Every minimal properly elliptic surface with good orbifold base,
odd first Betti number and Kodaira dimension 1 which does not have
singular fibres is the quotient of $\widetilde{SL(2,{\mathbb
R})}\times {\mathbb R}$ by a discrete subgroup \cite[Theorem
7.4]{Wall}. Also, every Inoue surface $S^+$ defined by means of a
real parameter $t$ is the quotient of the group $Sol_1^4$ by a
discrete subgroup \cite[Proposition 9.1]{Wall}. Left-invariant
para-hypercomplex structures descending to these types of elliptic
and Inoue surfaces have been constructed in \cite{AnSal, BV}.
Compatible metrics have been given in \cite{Iv-Zam} where it is
shown that the respective para-hyperhermitian structures are
locally, but not globally, conformally para-hyperk\"ahler, the
metric on the elliptic surfaces being locally conformally flat.
The flat para-hyperk\"ahler structures  on compact complex
surfaces have been described in \cite{Kamada99, Kamada02}; they
exist only on complex tori and primary Kodaira surfaces. Notice
that the structures constructed in Theorem~\ref{exist} on the
quaternionic Hopf surfaces and elliptic surfaces are  locally
conformally flat for $\varphi=const$.

\section{Locally conformally para-hyperherk\"ahler surfaces }

We have seen in the proof of Theorem~\ref{exist} that all surfaces
listed there admit locally conformally para-hyperk\"ahler
structures. The next result shows that these are the only compact
complex surfaces admitting such structures.

\begin{thm}\label {lcpk}
If a compact complex surface $(M,J)$ admits a locally conformally para-hyperk\"aher structure
$(g,I,S,T)$ with $I=J$ it is one of the following: a complex torus, a primary Kodaira surface, an
Inoue surface of type $S^{+}$, a properly elliptic surface of real type with odd first Betti number
or a Hopf surface of real type.
\end{thm}
{\bf Proof}. We have to show that some of the surfaces in
Theorem~\ref{Ch0} do not admit locally conformally
para-hyperk\"ahler structures. We first exclude the K3 surfaces.
All K3 surfaces are simply connected, so any locally conformally
para-hyperk\"aher structure on a K3 surface is globally
conformally para-hyperk\"aher and after a conformal change it
becomes para-hyperk\"ahler. However the K3 surfaces do not admit
such structures as proven by Kamada \cite{Kamada99, Kamada02}.

Let $M$ be an Inoue surface of type $S^-$  defined via $g_0(z,w)=(-z,\alpha w)$ and  $g_i$,
$i=1,2,3$, as in the proof of Theorem \ref{exist}. Assume that $M$ admits a locally conformally
para-hyperk\"aher structure and denote by $\Omega'$ and $\theta'$ the $(2,0)$-form and the Lee form
of this structure. For suitable $(p,q)\in\mathbb{Z}^2$ in the definition of Inoue surfaces, the
group generated by $g_0^2$,$g_1$,$g_2$,$g_3$ defines an Inoue surface $\widetilde{M}$ of type $S^+$
such that $M$ is the quotient of $\widetilde{M}$ by the fixed point free involution $\sigma$
determined by $g_0$ \cite[p. 279]{Inoue}. Denote by $\pi: \widetilde{M}\rightarrow M$ the
projection and  set $\widetilde{\Omega}=\pi^*\Omega'$ and $\widetilde{\theta}=\pi^*\theta'$. Then
$d\widetilde{\Omega}=\widetilde{\theta}\wedge\widetilde{\Omega}$ and $d\widetilde{\theta}=0$. Lift
the forms $\widetilde{\Omega}$ and $\widetilde{\theta}$ to the universal covering ${\mathbb
C}\times{\bf H}$ of $\widetilde M$ and denote the forms obtained by the same symbols. If
$\alpha_1,...,\alpha_4$ are the $G^+$-invariant $1$-forms defined by (\ref{s+}), then
$\Omega=(\alpha_1+i\alpha_2)\wedge (\alpha_3+i\alpha_4)$ is a nowhere-vanishing $G^+$-invariant
$(2,0)$-form on ${\mathbb C}\times{\bf H}$. Hence $\widetilde\Omega=f\Omega$ for a complex-valued
nowhere-vanishing smooth function. Then
$$
d\widetilde\Omega=df\wedge\Omega+fd\Omega=(\displaystyle{\frac{df}{f}}-\alpha_4)\wedge\widetilde\Omega
$$
Since ${\mathbb C}\times{\bf H}$ is simply connected, there is a smooth function $g$ such that $e^g=f$. If we set $\psi=Im\,g$, then $g=\ln|f|+i\psi$
and $\displaystyle{\frac{df}{f}}=d\ln|f|+id\psi$. We have $\partial\psi\wedge\widetilde\Omega=0$ since $\widetilde\Omega$ is of type $(2,0)$, thus
$$
\frac{df}{f}\wedge\widetilde\Omega=(d\ln|f|+i\overline{\partial}\psi-i\partial\psi)\wedge\widetilde\Omega=(d\ln|f|+d^c\psi)\wedge\widetilde\Omega
$$
Therefore
$$
\widetilde\theta\wedge\widetilde\Omega=(d\ln|f|+d^c\psi-\alpha_4)\wedge\widetilde\Omega.
$$
and it follows that
\begin{equation}\label{tet}
\widetilde\theta-(d\ln|f|+d^c\psi-\alpha_4)=0
\end{equation}
since the $1$-form on the left-hand side of (\ref {tet}) is
real-valued and $\widetilde\Omega$ is of type $(2,0)$. The
function $f$ is $G^+$-invariant since $\widetilde\Omega$ and
$\Omega$ are invariant. Hence the function $\ln|f|$ and the form
$d\psi$ are also $G^+$-invariant as well as $d^c\psi=-Jd\psi$, $J$
being the complex structure (but the function $\psi$ is not
necessarily invariant). Consider $\ln|f|$ and $d^c\psi$ on the
surface $\widetilde M$ and note that the form
$\widetilde\theta-d\ln|f|$ is closed on $\widetilde M$. It is
shown in \cite{Inoue} that the first Betti number of $\widetilde
M$ is equal to $1$. The form $\alpha_4$ considered on $\widetilde
M$ is closed and nowhere-vanishing, hence not exact. Therefore
there are a real constant $C$ and a real-valued smooth function
$\eta$ on $\widetilde M$ such that
\begin{equation}
\widetilde\theta-d\ln|f|=C\alpha_4+d\eta.
\end{equation}
Then, by (\ref{tet}),
\begin{equation}\label{dcpsi}
d^c\psi=(C+1)\alpha_4+d\eta
\end{equation}
Applying the operator $d^c$ to both sides, we obtain
$$
0=(C+1)d^c\alpha_4+d^cd\eta.
$$
It follows from (\ref{s+}) that $d^c\alpha_4=\alpha_3\wedge\alpha_4$. Hence
$$
d^cd\eta=-(C+1)\alpha_3\wedge\alpha_4.
$$
 Let $h$ be an Hermitian metric on $\widetilde M$ and denote by $\omega$ the fundamental $2$-form of the Hermitian manifold $(\widetilde M,h)$. Then
$$
h(d^cd\eta,\omega)=-(C+1)h(\alpha_3\wedge\alpha_4,\omega).
$$
Extend $J$ on $1$-forms by $J(\alpha)(X)=-\alpha(JX)$. Then $J\alpha_3=\alpha_4$ and we get
$$
\square_h\eta=-(C+1)|\alpha_4|^2_h
$$
where $\square_h$ is the complex Laplacian. The right-hand side of the latter identity has a constant sign, hence $\eta=const$ by the maximum principle. Therefore $C=-1$ and identity (\ref{dcpsi}) becomes $d^c\psi=0$ on $\widetilde M$, hence $d^c\psi=0$ on ${\mathbb C}\times{\bf H}$. Thefore $\psi=const=c$ and we have $\widetilde\Omega=|f|e^{ic}\Omega$. Consider the latter identity on $\widetilde M$. The form $\widetilde\Omega$ is $\sigma$-invariant,
while $\sigma^*(\Omega)=-\Omega$. It follows that $|f\circ\sigma|=-|f|$.  However $|f|$ is
positive everywhere, a contradiction.

Now we shall discuss the Inoue surfaces of type $S^0$ in a similar way. First recall their
construction. Let $A\in SL(3,\mathbb{Z})$ be a matrix with two complex eigenvalues $\alpha$ and
$\overline{\alpha}$, and a real eigenvalue $c>1$. Choose eigenvectors
$(\alpha_1,\alpha_2,\alpha_3)\in \mathbb{C}^3$ and $(c_1,c_2,c_3)\in\mathbb{R}^3$  corresponding to
$\alpha$ and $c$, respectively. Then the vectors $(\alpha_1,\alpha_2,\alpha_3),
(\overline{\alpha_1},\overline{\alpha_2},\overline{\alpha_3})$ and $(c_1,c_2,c_3)$ are
$\mathbb{C}$-linearly independent. Let $\Gamma$ be the group of automorphisms of ${\mathbb
C}\times{\bf H}$ generated by
$$g_o:(z,w)\to (\alpha z,cw), \quad g_i:(z,w)\to (z+\alpha_i,w+c_i), \>i=1,2,3 .$$
Then $S=({\mathbb C}\times{\bf H})/\Gamma$ is an Inoue surface of type $S^0$ \cite{Inoue}.

Set $w=u+iv$, $a=\ln|\alpha|$, $b=-Arg\alpha$, $0<Arg\alpha\leq2\pi$, and $t=\displaystyle{\frac{\ln v}{\ln
c}}$. Define real vector fields $E_1,...,E_4$ on ${\mathbb C}\times{\bf H}$ by
\begin{equation}\label{1.1}
E_1-iE_2=2\alpha^t\frac{\partial}{\partial z},\hspace{.1in}
E_3-iE_4=2v\ln c\frac{\partial}{\partial w}.
\end{equation}
 These vector fields are  $\Gamma$-invariant, hence they define $(1,0)-$ vector fields on
 $S$ which we denote by the same symbols. Thus, if $J$ is the complex structure of $S$, then  $JE_1=E_2$ and $JE_3=E_4$. It is easy to check that
\begin{equation}\label{1.2}
[E_4,E_1]=aE_1-bE_2,\quad [E_4,E_2]=bE_1+aE_2, \quad [E_4,E_3]=-2aE_3
\end{equation}
and all other brackets vanish. Denote by $\alpha_1,...,\alpha_4$ the
dual frame of $E_1,...,E_4$. Clearly, the $1$-forms $\alpha_i$ are $\Gamma$-invariant. Moreover, (\ref {1.2}) implies that

\begin{equation}\label{1.3}
\begin{array}{c}
d\alpha_1=a\alpha_1\wedge\alpha_4+b\alpha_2\wedge\alpha_4,\quad
d\alpha_2=-b\alpha_1\wedge\alpha_4+a\alpha_2\wedge\alpha_4,\\[6pt]
d\alpha_3=-2a\alpha_3\wedge\alpha_4, \quad d\alpha_4=0.
\end{array}
\end{equation}

Suppose that $S$ admits a locally conformally para-hyperk\"aher
structure. Let $\Omega'$ and $\theta'$ be the $(2,0)$-form and the
Lie form of this structure. Then $d\Omega'=\theta'\wedge\Omega'$
and $d\theta'=0$. We lift $\Omega'$ and $\theta'$ to the universal
covering  ${\mathbb C}\times{\bf H}$  of $S$ and denote the lifts
by the same symbols. Set $\Omega=(\alpha_1+i\alpha_2)\wedge
(\alpha_3+i\alpha_4)$. Then $\Omega'=f\Omega$ for a complex-valued
nowhere-vanishing smooth function $f$ on ${\mathbb C}\times{\bf
H}$. We have $d\Omega'=df\wedge\Omega+fd\Omega$ where, in view of
(\ref{1.3}), $d\Omega=(-b\alpha_3+a\alpha_4)\wedge\Omega$. Thus
\begin{equation}\label{d-om}
d\Omega'=(\frac{df}{f}-b\alpha_3+a\alpha_4)\wedge\Omega'.
\end{equation}
Since ${\mathbb C}\times{\bf H}$ is simply connected, there is a smooth function $\psi$ such that
$f=|f|e^{i\psi}$. We have $\partial\psi\wedge\Omega'=0$ since $\Omega'$ is of type $(2,0)$. Then
$$
\frac{df}{f}\wedge\Omega'=(d\ln|f|+i\overline{\partial}\psi-i\partial\psi)\wedge\Omega'=(d\ln|f|+d^c\psi)\wedge\Omega'
$$
and it follows from (\ref{d-om}) that
$$
\theta'\wedge\Omega'=(d\ln|f|+d^c\psi-b\alpha_3+a\alpha_4)\wedge\Omega'.
$$
This implies
\begin{equation}\label{thet}
\theta'-(d\ln|f|+d^c\psi-b\alpha_3+a\alpha_4)=0.
\end{equation}
The form $\theta'-d\ln|f|$ is $\Gamma$-invariant, so we can consider it on the surface $S$. According to \cite{Inoue}, $b_1(S)=1$. The form $\alpha_4$
considered on $S$ is closed and is not exact. Therefore there are a real constant $C$ and a real-valued smooth function $\eta$ on $S$ such that
\begin{equation}\label{onS}
\theta'-d\ln|f|=C\alpha_4+d\eta.
\end{equation}
Identities (\ref{thet}) and (\ref{onS}) imply that
\begin{equation}\label{dc}
d^c\psi=b\alpha_3+(C-a)\alpha_4+d\eta.
\end{equation}
Applying the operator $d^c$ to both sides, we obtain
$$
0=bd^c\alpha_3+(C-a)d^c\alpha_4+d^cd\eta.
$$
It follows from (\ref{1.1}) and (\ref{1.2}) that $d^c\alpha_3=0$ and $d^c\alpha_4=-2a\alpha_3\wedge\alpha_4$. Thus
$$
d^cd\eta=2a(C-a)\alpha_3\wedge\alpha_4.
$$
Take an Hermitian metric $h$ on $S$ (for example that for which $E_1,...,E_4$ is an orthonormal frame) and let $\omega$ be the fundamental $2$-form of the Hermitian manifold $(S,h)$. Then
$$
h(d^cd\eta,\omega)=2a(C-a)h(\alpha_3\wedge\alpha_4,\omega).
$$
and, since $J\alpha_3=\alpha_4$, and we get
$$
\square_h\eta=2a(C-a)|\alpha_4|^2_h.
$$
The latter identity implies $\eta=const$. Therefore $C=a$ and identity (\ref{dc}) takes the form $d^c\psi=b\alpha_3$.
Since $d^c\psi=Jd\psi$, we get $d\psi=-b\alpha_4$. According to (\ref{1.1}), $\alpha_4=\displaystyle{\frac{dv}{v\ln c}}$, hence
$\displaystyle{d\psi=d(b\frac{\ln v}{\ln c})}$.
Therefor there is a constant $C'$ such that
$$
\exp(i\psi)=C'\exp(ib\frac{\ln v}{\ln c}).
$$
Since $exp(i\psi)=\displaystyle{\frac{f}{|f|}}$, the function on the right-hand side of the latter identity is $\Gamma$-invariant. In particular, this function is invariant under the transformation $g_0$, hence
$$
exp(ib\frac{\ln(cv)}{\ln c})=\exp(ib\frac{\ln v}{\ln c}).
$$
This gives $\exp(ib)=1$, hence $b=2k\pi$ for an integer $k$. But this means that the eigenvalue
$\alpha$ of the matix $A$ is a real number, a contradiction.    {\it Q.E. D.}

%Finally, we note that, as we have seen in the proof of Theorem~\ref{exist}, every complex torus,
%primary Kodaira surface, Inoue surface of type $S^{+}$, as well as every special type minimal
%properly elliptic surface with odd first Betti number and every quaternionic primary Hopf surface
%does admit a locally conformally para-hyperk\"ahkler structure.

 %The surfaces of type $S^-$ are factors of a surface of type $S^+$ by a holomorphic involution $\sigma$. They admit a non-vanishing $(2,0)$-form %$\Omega$, such that $\Omega$ pulls-back to an invariant $(2,0)$-form $\widetilde{\Omega}$ on the covering surface, such that %$\sigma^*(\widetilde{\Omega})=-\widetilde{\Omega}$ - see \cite{DGMY1}. But such form can not be part of a locally conformally para-hyperk\"ahler %structure from Lemma~\ref{real}. Assume that $\Omega_2+i\Omega_3 = fe^{ig}\Omega = \Omega'$ is a $(2,0)$-form of a locally conformally %para-hyperk\"ahler structure on a given surface of type $S^-$. Then the pullback $\widetilde{\Omega'}$ of $\Omega'$ is a
%locally conformally para-hyperk\"ahler structure on the covering surface and $\widetilde{\Omega'} = \widetilde{f}e^{\widetilde{g}}\widetilde{\Omega}$ %for some functions $\widetilde{f}>0$ and $\widetilde{g}$ which are pullbacks of $f$ and $g$ respectively. But then $\sigma^*(fe^{ig}) = %\widetilde{f}e^{i\widetilde{g}}$ and $\sigma^*(\Omega') = -\Omega'$. A contradiction.

\section{Para-hypercomplex surfaces}

Now we shall use Theorem \ref{Ch0} to provide a list of the compact
complex surfaces that could admit a para-hypercomplex structure.

\begin{thm}Let $(M,J)$ be a compact complex surface admitting a
para-hypercomplex structure. Then it is one of the surfaces listed
in Theorem~\ref{Ch0}, a hyperelliptic surface, a secondary Kodaira
surface, or an Enriques surface.
\end{thm}
{\bf Proof}. Every complex surface with a para-hypercomplex
structure which does not admit a para-hyperhermitian structure has
a double cover which admits a para-hyperhermitian structure
compatible with the pull-back of the para-hypercomplex structure
(Proposition~ \ref{cover}). Then it follows from the list of
possible para-hyperhermitian surfaces in Theorem~\ref{Ch0} that we
have to consider only that admitting holomorphic involutions and
to identify the corresponding quotient surfaces. It is well known
that a smooth quotient of a torus, a K3 surface, or a primary
Kodaira surface is, respectively, a hyperelliptic surface, an
Enriques surface, or a secondary Kodaira surface. Also the
quotient of an Inuoe surface with $b_2=0$ by a holomorphic
involution is a surface of the same type. Note also that such a
quotient of a Hopf surface is a Hopf surface since it has the same
universal cover.

Let $\pi: M\rightarrow C$ be a properly elliptic surface with odd $b_1$ and base $C$ of genus
$g>1$. Then $M$ does not have singular fibers and multisections. If $M$ admits an involution
$\tau$, then $\tau$ transforms a fixed fiber $E$ into a curve $C'$.  Any curve in $M$ which is not
a fiber should project onto the whole base and hence should be a multisection, a contradiction. So
the projection $\pi(C')$ is a point, hence the image $\tau(E)=C'$ is contained in a fiber. Since
all fibers are irreducible smooth elliptic curves, $C'$ is again a fiber, possibly with different
multiplicity. Then $\tau$ induces an involution $\tau'$ of the base $C$ and $M/\tau$ is
elliptically fibred over $C/\tau'$ without singular fibers. It should have vanishing real first
Chern class. Then it is either a properly elliptic surface, a Hopf surface or a Kodaira surface
depending on the genus of $C/\tau'$. {\it Q.E.D.}

\smallskip

In \cite{DGMY} we have shown that every Inoue surface of type
$S^{-}$ has a para-hypercomplex structure which does not admit a
compatible para-hyperhermitian metric. Here we construct such a
structure on a hyperelliptic surface.

\smallskip

\noindent {\bf Example}. Let $T^2=\mathbb{C}/<1,i>$ be the complex torus with lattice generated by
$1$ and $i$. Denote by $\varphi$ the holomorphic involution of $T^2\times T^2$ defined by $(z,w)\to
(z+\frac{i}{2},-w)$. Then the quotient $M$ of $N=T^2\times T^2$ by the group generated by
${\varphi}$ is a hyperelliptic surfaces. Let $z=x+iy$,  resp. $w=u+iv$ be the local coordinate on
the first, resp. the second factor of $N$ induced by the standard complex coordinate of
$\mathbb{C}$. Then
$$
K^{+}=span\{\frac{\partial}{\partial x},\frac{\partial}{\partial
u}\}, \quad K^{-}=span\{\frac{\partial}{\partial
y},\frac{\partial}{\partial v}\}
$$
are $\varphi$-invariant and involutive subbundles of the tangent bundle $TN$. Define an isomorphism
$S$ of $TN$ setting $S=+Id$ on $K^{+}$ and $S=-Id$ on $K^{-}$. Let $I$ be the complex structure of
$N$ and set $T=IS$. In this way we obtain a para-hypercomplex structure on $N$ which descends to
$M=N/<\varphi>$. This structure does not admit a compatible metric since otherwise, as we have
seen, the canonical bundle of the hyperelliptic surface $M$ would be topologically trivial, a
contradiction with Theorem 8.

\section{Nonexistence of para-hyperhermitian structures on Inoue surfaces of type $S^0$}

It is well-known \cite{Hasegawa} that every Inoue surface of type
$S^0$ is a solvmanifold, i.e. the quotient of a solvable Lie group
by a cocompact subgroup. Note also that the 4-dimensional solvable
Lie algebras admitting para-hypercomplex structures have been
classified in \cite{BV}. This together with the identities (\ref
{1.2}) above implies that the Inoue surfaces of type $S^0$ do not
admit para-hyperhermitian structures induced by left invariant
ones. In this section we shall slightly generalize this
observation by showing that these surfaces have no
para-hyperhermitian structures whose $(2,0)$-forms are defined by
left invariant 2-forms. This leads to the natural conjecture that
the Inoue surfaces of type $S^0$ do not admit para-hyperhermitian
structures at all.

%In what follows, we shall use facts and formulas stated in the proof of Theorem~\ref{lcpk}.

\begin{thm}\label{Inoue} Let $S$ be an Inoue surface of type $S^0$
which is a quotient of a solvable Lie group $G$. Then $S$ has no
para-hyperhermitian structure whose $(2,0)$-form is defined by a
left invariant 2-form on $G$.
\end{thm}
{\bf Proof.}  We shall use the notation introduced in the proof of Theorem~\ref{lcpk}.
To prove the theorem, we have to consider the class of para-hyperhermitian structures on $S$ whose
$(2,0)$-form $\Omega_2+i\Omega_3$ is given (up to a constant) by

$$\Omega_2+i\Omega_3=
(\alpha_1+i\alpha_2)\wedge(\alpha_3+i\alpha_4),$$ i.e.
$$\Omega_2=\alpha_1\wedge\alpha_3-\alpha_2\wedge\alpha_4,\quad
\Omega_3=\alpha_1\wedge\alpha_4+\alpha_2\wedge\alpha_3.$$ Since
$\Omega_1$ is a real $(1,1)$-form with respect to $J$ it has the
form
$$\Omega_1=p\alpha_1\wedge\alpha_2+q(\alpha_1\wedge\alpha_3+\alpha_2\wedge\alpha_4)+
r(\alpha_1\wedge\alpha_4 -
\alpha_2\wedge\alpha_3)+s\alpha_3\wedge\alpha_4,$$ where $p,q,r,s$
are smooth functions on $S$. Further on, we shall  consider the
smooth functions on $S$ as smooth $\Gamma-$invariant functions on
${\mathbb C}\times{\bf H}$. Denote by $\theta$ the Lie form of the para-hyperhermitian
structure on $S$ and set $f=q+ir$. Then a direct but long
computation using (\ref {1.1}) and (\ref {1.3}) shows that the
identities in Proposition \ref{phe} are satisfied if and only if
$\theta=-b\alpha_3+a\alpha_4$, $p=0$, $|f|^2=1$ and

$$\frac{\partial f}{\partial z}=0,\quad
d\overline{f}(c^t\frac{\partial}{\partial
w})=-2b\overline{f}-ids(\alpha^t\frac{\partial}{\partial z}).$$ Now
differentiating the identity $|f|^2=1$ with respect to
$\overline{z}$ we get $\ds\frac{\partial f}{\partial\overline z}=0.$
Hence the function $f$ depends only on $u$ and $v$ and satisfies the
identity
\begin{equation}\label{iden}
\ln c\frac{\partial f}{\partial\overline
w}=-bf-i\overline{ds(\alpha^t\frac{\partial}{\partial z})}.
\end{equation}

Next we shall need the following
\begin{lm}Let $F$ be a continuous function on ${\mathbb C}\times{\bf H}$ which is
invariant under the action of $\Gamma$ and depends only on $u$ and
$v$. Then $F$ depends only on $v$.
\end{lm}
{\bf Proof.} The invariance of $F$ implies that
\begin{equation}
F(u+xc_1+yc_2+zc_3,v)=F(u,v)
\end{equation}
for arbitrary $x,y,z\in\mathbb{Z}.$ Note that at least two of the numbers $c_1,c_2,c_3$ are nonzero
and since the eigenvalue $c$ is irrational we may assume that the ratio
$\displaystyle{\frac{c_1}{c_2}}$ is irrational too. Then the Kronecker lemma implies that the set
$\{xc_1+yc_2|~x,y\in \mathbb{Z}\}$ is dense in $\mathbb{R}$ and by continuity we get that
$F(u+w,v)=F(u,v)$ for arbitrary $u,v,w\in \mathbb{R}$. Hence $F$ depends only on $v$.  {\it Q. E.
D.}

The lemma above implies that the function $f$ depends only on $v$ and it follows from (\ref{iden})
that the same is also true for $A=\displaystyle{\frac{\partial s}{\partial z}}$  . Since $s$ is a
real function we have $\displaystyle{\frac{\partial s}{\partial \overline{z}}=\overline{A}}$ and
therefore
\begin{equation}\label{eq1}
s=zA+\overline{z}\overline{A}+\gamma(u,v).
\end{equation}
The invariance of $s$ implies that
$$s(z+x\alpha_1+y\alpha_2+t\alpha_3,\overline{z}+x\overline{\alpha_1}
+y\overline{\alpha_2}+t\overline{\alpha_3},u+xc_1+yc_2+tc_3,v)=
s(z,\overline{z},u,v)$$ for all $x,y,t\in\mathbb{Z}$ and it
follows from (\ref{eq1}) that
\begin{equation}\label{eq2}
xA_1+yA_2+tA_3=\gamma(u,v)-\gamma(u+xc_1+yc_2+tc_3,v),
\end{equation}
where $A_i=\alpha_iA+\overline{\alpha_i}\overline{A}$, $i=1,2,3.$

We shall show now that $A=0$. To do this we consider two cases.

{\it Case 1.} $c_1c_2c_3\neq0.$

Take sequences $\{x_n\}$ and $\{y_n\}$ of integers such that
$\displaystyle{x_n\frac{c_1}{c_2}+y_n}$ tends to a rational number $\omega$. Then
setting $x=x_n$, $y=y_n$, $t=0$ and $u=0$ in (\ref{eq2}) gives

$$x_nA_1+y_nA_2=\gamma(0,v)-\gamma(x_nc_1+y_nc_2,v).$$
Hence $x_nA_1+y_nA_2$ tends to $B=\gamma(0,v)-\gamma(c_2\omega,v)$. We may assume without loss of
generality that $A_2>0$. Then there exists $N$ such that for every $n>N$ we have
$$|x_n\frac{c_1}{c_2}+y_n-\omega|<1.$$ On the other hand there exists $M$ such that for every $n>M$
we have
$$|x_nA_1+y_nA_2-B|<1.$$ The last two inequalities
imply that for $n>\max(M,N)$ we have
$$B+1-A_2(1+\omega)<x_n(A_1-\frac{A_2c_1}{c_2})<B+1-A_2(1-\omega).$$
Suppose that $\displaystyle{A_1-\frac{A_2c_1}{c_2}\neq0}$. Since $x_n$ are
integers this inequality implies  that $x_n$ take a
finite number of values and the same is also true for $y_n$. But
the sequence $\displaystyle{\{x_n\frac{c_1}{c_2}+y_n\}}$ is convergent and
therefore its limit is equal to a term of it, a contradiction
since the number $\displaystyle{\frac{c_1}{c_2}}$ is irrational. Thus
$\displaystyle{A_1(\omega)=\frac{c_1}{c_2}A_2(\omega)}$ for every
$\omega\in\mathbb{Q}$, hence for every $\omega\in\mathbb{R}$.
The same reasoning shows that the vectors $(A_1,A_2,A_3)$ and
$(c_1,c_2,c_3)$ are collinear and therefore the vectors
$(A\alpha_1,A\alpha_2,A\alpha_3),
(\overline{A}\overline{\alpha_1},\overline{A}\overline{\alpha_2},\overline{A}\overline{\alpha_3})$
and $(c_1,c_2,c_3)$ are $\mathbb{C}$-linearly dependent. Hence
$A=0$.

{\it Case 2.} $c_1c_2c_3=0.$

 We may assume that $c_3=0$, $c_1c_2\neq0$. Then, applying (\ref{eq2}) for $x=y=0$, $t=1$, we
get $A_3=0$. Since $c_2\neq0$ the same reasoning as in Case 1
implies that the vectors $(A_1,A_2,A_3)$ and $(c_1,c_2,c_3)$ are
collinear and we get again that $A=0$.

Now the equation (\ref{iden}) takes the form $$v\ln c\frac{\partial f}{\partial v}=2ibf.$$ Since
$|f|=1$ and $f$ depends only on $v$, it follows that $f=e^{ig}$, where $g$ is a smooth real-valued
function on $\mathbb{R}_+$. Then the latter equation takes the form
$$\ds{\frac{\partial g}{\partial v}=\frac{2b}{v\ln c}}$$ which shows
that
$$g=\frac{2b\ln v}{\ln c} +g_0,$$ where $g_0$ is a constant.
Therefore $$f(v)=f_0\exp(\ds\frac{2ib\ln v}{\ln c}),$$ where $f_0$
is a constant with $|f_0|=1$. Now the invariance of $f$ under
$\Gamma$ implies that $b=k\pi, k\in \mathbb Z$. Hence $\alpha$ is a
real number, a contradiction. {\it Q.E.D.}

 \vskip 20pt

%\newpage

\noindent Johann Davidov

\noindent Institute of Mathematics and Informatics

\noindent Bulgarian Academy of Sciences

\noindent 1113 Sofia, Bulgaria

\noindent

and

\noindent "L.Karavelov" Civil Engineering Higher School

\noindent 175 Suhodolska Str.

\noindent 1373 Sofia, Bulgaria

\noindent jtd@math.bas.bg

\medskip

\noindent Oleg Mushkarov

\noindent Institute of Mathematics and Informatics

\noindent Bulgarian Academy of Sciences

\noindent 1113 Sofia, Bulgaria

\noindent muskarov@math.bas.bg

\medskip

\noindent Gueo Grantcharov and Miroslav Yotov

\noindent Department of Mathematics

\noindent Florida International University

\noindent Miami, FL 33199

\noindent grantchg@fiu.edu, yotovm@fiu.edu

\end{document}